\documentstyle{article}

\title{Solvability of the boundary value problem\\ for some
differential--difference equations}
\author{Pavel~Gurevich}
\date{}

\begin{document}
\maketitle
\begin{abstract}
Solvability and smoothness of generalized solutions to boundary
value problems for not self-adjoint differential-difference
equations are studied. Necessary and sufficient conditions of
Fredholmian solvability (with index zero) are established.
Smoothness of generalized solutions is considered in terms of
index of the corresponding differential--difference operator.
\end{abstract}

\begin{center}
\large\bf Inroduction
\end{center}

The problems of solvability and smoothness of generalized
solutions to boundary value problems for differential--difference
equations on a finite interval $(0,d)$ in not self-adjoint case
were considered in [3]. The interest to these problems was arised
by their numerous applications as well as by a number of quite new
properties they possess. For instance, the smoothness of
generalized solutions to such problems may fail inside the
interval $(0,d)$ even in the case of infinitely differentiable
right hand side of the equation and remains only in some
subintervals. In [3] necessary and sufficient conditions of
Fredholmian solvability and smoothness of solutions to such
problems on the whole interval were established in the case of
non--integer $d$. In the case of integer $d$ only sufficient
conditions were obtained. The problem of obtaining necessary and
sufficient conditions was formulated in [3] as an unsolved one.
This paper is dedicated to the solution of this problem.

In \S1, the properties of difference operators in Sobolev spaces
are considered. In \S2, the necessary and sufficient conditions of
Fredholmian solvability (with index zero) of a boundary value
problem for a differential--difference equations are established.
In \S3, the smoothness of the generalized solutions is considered
in terms of the index of the corresponding
differential--difference operator.

The author is grateful to professor A.\,L.\,Skubachevski\v{\i} for
constant attention to this work.

\begin{center}
\large\bf\S1. Difference operators in the spaces $L_2({\bf R}),$ $L_2(0,N+1),$\\
and in the Sobolev spaces $W^k(0,N+1)$
\end{center}
We consider the {\it difference} operator $R:L_2({\bf R})\to L_2({\bf R})$
defined by the formula
$$
  (Rv)(t)=\sum\limits_{j=-N}^N b_j y(t+j).\eqno(1)
$$
Here $b_j$ are real numbers, $N$ is a natural number.

We introduce the operators
$$
  \begin{array}{l}
    I_Q:L_2(0,N+1)\to L_2({\bf R}),\quad P_Q:L_2({\bf R})\to L_2(0,N+1),\\
    R_Q:L_2(0,N+1)\to L_2(0,N+1)
  \end{array}
$$
by the formulas
$$  
  \array{l}   
    (I_Qv)(t)=\left\{\begin{array}{ll} v(t) & (t\in(0,N+1)),\\
	           	              0    & (t\not\in(0,N+1));
              \end{array}\right.\qquad  (P_Qv)(t)=v(t)\quad (t\in(0,N+1));\\
    R_Q=P_QRI_Q.
  \endarray\eqno(2)
$$
Here $Q=(0,N+1).$

We denote $Q_s=(s-1,s)\quad (s=1,\dots,N+1).$

We introduce an isomorphism of the Hilbert spaces
$$
  U:L_2(\cup_sQ_s)\to L_2^{N+1}(Q_1)
$$
by the formula
$$
  (Uv)_k(t)=v(t+k-1)\quad (t\in Q_1,\; k=1,\dots,N+1),\eqno(3)
$$
where $L_2^{N+1}(Q_1)=\prod_{k=1}^{N+1} L_2(Q_1).$

Let $R_1$ be the matrix of order $(N+1)\times (N+1)$ with the elements
$
  r_{ik}=b_{k-i}\quad (i,k=1,\dots,N+1).
$
Let $R_2$ be the matrix of order $N\times N$ obtained from $R_1$ by deleting
the last column and the last row.
We denote also by $B_{ik}$ the cofactor of the element $r_{ik}$ of the
matrix $R_1.$

Consider the operator $R_{Q1}:L_2^{N+1}(Q_1)\to L_2^{N+1}(Q_1)$ defined by the
formula $R_{Q1}=UR_QU^{-1}.$

Now we shall formulate the next four Lemmas (proofs are given in [3],
Chapter I, Section 2).
\bigskip

{\bf Lemma 1.} {\sl The operator $R_{Q1}$ is the operator of
multiplication by the matrix $R_1.$}
\bigskip

{\bf Lemma 2.} {\sl The spectrum of the operator $R_Q$ coincides with the
spectrum of the matrix $R_1.$}
\bigskip

{\bf Lemma 3.} {\sl The operator $R_Q$ maps continuously $\mathaccent23W^k(0,N+1)$
into $W^k(0,N+1)$ and, for all $v\in\mathaccent23W^k(0,N+1),$
$$
  (R_Qv)^{(j)}=R_Qv^{(j)}\quad (j\le k).\eqno(4)
$$
}
\bigskip

{\bf Lemma 4.} {\sl Let $\det R_1\ne0$ and let $R_Qv\in W^k(Q_i)$ for
$i=1,\dots,N+1.$
Then $v\in W^k(Q_j)\;(j=1,\dots,N+1)$ and
$$
  \|v\|_{W^k(Q_j)}\le c\sum\limits_{i=1}^{N+1} \|R_Qv\|_{W^k(Q_i)},
$$
where $c>0$ doesn't depend on $v.$ }
\bigskip

Let us denote by $W^k_\gamma(0,N+1)$ the subspace of functions from $W^k(0,N+1)$
satisfying conditions
$$
  u^{(\mu)}(N+1)=\sum\limits_{1\le i\le N+1,\, i\ne m+1}\gamma_{1i}u^{(\mu)}(i-1),
								       \eqno(5)
$$
$$
  u^{(\mu)}(m)=\sum\limits_{1\le i\le N,\,i\ne m}\gamma_{2i}u^{(\mu)}(i),\eqno(6)
$$
where $m$ is a fixed point from the set $\{1,\dots,N\},$
$\gamma_{1i}(i=1,\dots,N+1,\, i\ne m+1),$ 
$\gamma_{2i}(i=1,\dots,N,\,i\ne m)$ are real numbers; 
$\mu=0,\dots,k-1;\; k\ge1.$

Hereinafter, we shall assume that $\det R_1\ne0,\;\det R_2=0$ as the other cases
have been studied in [3], Chapter I.

\bigskip
{\bf Theorem 1.} {\sl There exist real numbers 
$\gamma_{1i}(i=1,\dots,N+1,\, i\ne m+1),$ 
$\gamma_{2i}(i=1,\dots,N,\,i\ne m)$ 
such that the operator $R_Q$ maps $\mathaccent23W^k(0,N+1)$ onto $W^k_\gamma(0,N+1)$
continuously and in a one--to--one manner.}
\bigskip

{\it Proof.} 1. At first we proof that there exist
$\gamma_{1i}(i=1,\dots,N+1,\, i\ne m+1),$ 
$\gamma_{2i}(i=1,\dots,N,\,i\ne m)$
such that $R_Q(\mathaccent23W^k(0,N+1))\subset W^k_\gamma(0,N+1).$

We denote by $R_1^1(R_1^2)$ the matrix, obtained from $R_1$ by deleting
the first (the last) column. Denote by $e_i\,(g_i)$ the $i$-th row of the
matrix $R_1^1(R_1^2).$ 

The condition $\det R_2=0$ implies that $g_1,\dots,g_N$ are linearly dependent.
Hence there exists a point $m$ from the set $\{1,\dots,N\}$ such that the row
$g_m$ is a linear combination of the other ones
$$
  g_m=\sum\limits_{1\le i\le N,\,i\ne m}\gamma_{2i}g_i,\eqno(7)
$$
where $\gamma_{2i}\,(i=1,\dots,N,\,i\ne m)$ are real numbers.

It is easy to see that $e_{i+1}=g_i\; (i=1,\dots,N).$ Therefore, using (7),
we get
$$
  e_{m+1}=\sum\limits_{1\le i\le N,\,i\ne m}\gamma_{2i}e_{i+1},\eqno(8)
$$
i.e.,
$$
e_{m+1}=\sum\limits_{2\le i\le N+1,\,i\ne m+1}\gamma_{2,i-1}e_i.\eqno(9)
$$ 

From the non-singularity of the matrix $R_1$ it follows that the rows
$e_i\;(i=1,\dots\,N+1,\,i\ne m+1)$ form the basis in ${{\bf R}}^N$ and the rows
$g_j\;(j=1,\dots\,N+1,\,j\ne m)$ do the same.

By Lemma 3, $R_Q(\mathaccent23W^k(0,N+1))\subset W^k(0,N+1).$ Thus (3), (7) and Lemma 1
implies that, for $v\in\mathaccent23W^k(0,N+1)$ and $\mu=0,\dots,k-1,$
$$
\begin{array}{ll}  
(R_Qv)^{(\mu)}(m)&=(UR_Qv)^{(\mu)}_m(1)                             \\
{ } &=(R_1Uv^{(\mu)})_m(1)=\sum\limits_{1\le i\le N,\,i\ne m} 
   					\gamma_{2i}(R_1Uv^{(\mu)})_i(1) \\
{ } &=\sum\limits_{1\le i\le N,\,i\ne m}\gamma_{2i}(R_Qv)^{(\mu)}(i) \qquad 
  (\mu=0,\dots,k-1).
\end{array}\eqno(10)
$$

Further,
$$
\begin{array}{ll}
  (R_Qv)^{(\mu)}(N+1)&=(UR_Qv)^{(\mu)}_{N+1}(1)                      \\
{ }  &=(R_1Uv^{(\mu)})_{N+1}(1)=\sum\limits_{s=1}^N r_{N+1,s}(Uv^{(\mu)})_s(1) \\
{ }  &=\sum\limits_{s=1}^N r_{N+1,s}(Uv^{(\mu)})_{s+1}(0)=
    \sum\limits_{s=2}^{N+1}r_{N+1,s-1}(Uv^{(\mu)})_s(0).
\end{array}\eqno(11)
$$

And, in the same way,
$$
\begin{array}{ll}
  (R_Qv)^{(\mu)}(i-1)&=(UR_Qv)^{(\mu)}_i(0)=(R_1Uv^{(\mu)})_i(0) \\
{ }  &=\sum\limits_{s=2}^{N+1}r_{is}(Uv^{(\mu)})_s(0)\quad (i=1,\dots,N+1). 
\end{array}\eqno(12)
$$

Since the rows $e_i\, (i=1,\dots,N+1;\,i\ne m+1)$ form the basis in 
${{\bf R}}^N$, it follows that 
$$
  g_{N+1}=\sum\limits_{1\le i\le N+1,\,i\ne m+1}\gamma_{1i}e_i,
$$
i.e.,
$$
  (r_{N+1,1},\dots,r_{N+1,N})=
  \sum\limits_{1\le i\le N+1,\,i\ne m+1}\gamma_{1i}(r_{i2},\dots,r_{i,N+1}).
								     \eqno(13)
$$

Now, using (11), (12), (13), we get
$$
  (R_Qv)^{(\mu)}(N+1)=
  \sum\limits_{1\le i\le N+1,\,i\ne m+1}\gamma_{1i}(R_Qv)^{(\mu)}(i-1)
  \quad (\mu=0,\dots,k-1).\eqno(14)
$$

Therefore, by virtue (10) and (14), $R_Q(\mathaccent23W^k(0,N+1))\subset W^k_\gamma(0,N+1).$

2. Now let us prove the inverse inclusion 
$$
W^k_\gamma(0,N+1)\subset R_Q(\mathaccent23W^k(0,N+1)).
$$

Suppose $u\in W^k_\gamma(0,N+1).$ By virtue of Lemma~2, the operator 
$R_Q:L_2(0,N+1)\to L_2(0,N+1)$ has a bounded inverse 
$R_Q^{-1}:L_2(0,N+1)\to L_2(0,N+1).$ We shall show that 
$v=R_Q^{-1}u\in\mathaccent23W^k(0,N+1).$

By virtue of Lemma 4, $v\in W(Q_s)\,(s=1,\dots,N+1).$ Therefore, to prove
this theorem, it is sufficient to prove that 
$$
  (Uv)_s^{(\mu)}(1-0)=(Uv)_{s+1}^{(\mu)}(0+0)\quad (s=1,\dots,N);
  \qquad (Uv)_1^{(\mu)}(0+0)=(Uv)_{N+1}^{(\mu)}(1-0)=0.
$$

Denote
$$
  \begin{array}{l}
   \varphi^\mu_s=(Uv)_{s+1}^{(\mu)}(0+0) \quad(s=0,\dots,N;\,\mu=0,\dots,k-1);\\
   \psi^\mu_j=(Uv)_j^{(\mu)}(1-0)    \quad(j=1,\dots,N+1;\,\mu=0,\dots,k-1).
  \end{array}
$$

Since $R_Qv\in W^k(0,N+1),$ we have
$$
  (R_Qv)^{(\mu)}|_{t=i-0}=(R_Qv)^{(\mu)}|_{t=i+0} \quad
  (i=1,\dots,N;\, \mu=0,\dots,k-1).
$$
Thus, for every $\mu=0,\dots\,k-1,$ the functions $\varphi^\mu_s,\; \psi^\mu_j$ 
satisfy the following
conditions
$$
  \sum\limits_{s=1}^{N+1}r_{i+1,s}\varphi^\mu_{s-1}=
  \sum\limits_{s=1}^{N+1}r_{is}\psi^\mu_s \quad (i=1,\dots,N). \eqno(15)
$$
Moreover, the function $R_Qv$ satisfies conditions (10), which can be
rewritten in the form
$$
  \sum\limits_{s=1}^{N+1}r_{ms}\psi^\mu_s=
  \sum\limits_{1\le i\le N,\,i\ne m} \gamma_{2i}\sum\limits_{s=1}^{N+1}r_{is}\psi^\mu_s, 
								\eqno(16)
$$
$$
  \begin{array}{ll}
           \sum\limits_{s=1}^{N+1}r_{m+1,s}\varphi^\mu_{s-1}&=
           \sum\limits_{1\le i\le N,\,i\ne m}\gamma_{2i}
				\sum\limits_{s=1}^{N+1}r_{i+1,s}\varphi^\mu_{s-1} \\
{ }        &=\sum\limits_{2\le i\le N+1,\,i\ne m+1}\gamma_{2,i-1}
           \sum\limits_{s=1}^{N+1}r_{is}\varphi^\mu_{s-1}.
  \end{array}\eqno(17)
$$

From conditions (16), (17) and (7), (9), we obtain
$$
  \begin{array}{l}
    \left(r_{m,N+1}-\sum\limits_{1\le i\le N,\,i\ne m}\gamma_{2i}r_{i,N+1}\right)\, 
							\psi^\mu_{N+1}=0,\\
    \left(r_{m+1,1}-\sum\limits_{2\le i\le N+1,\,i\ne m+1}
                   \gamma_{2,i-1}r_{i1}\right)\, \varphi^\mu_0=0.  
  \end{array}
$$

The factor preceding $\psi^\mu_{N+1}\, (\varphi^\mu_0)$ is non--zero. 
Otherwise, we have $\det R_1=0,$ which contradicts the conditions of the
theorem.
Hence $\psi^\mu_{N+1}=\varphi^\mu_0=0.$

Thus system (15) will have the form
$$
  \sum\limits_{s=1}^Nr_{i+1,s+1}\varphi^\mu_s=
  \sum\limits_{s=1}^N r_{is}\psi^\mu_s \quad (i=1,\dots,N).
$$
Since $r_{i+1,s+1}=r_{is}$ and the $m$-th row of this system is a linear
combination of the other ones, this system will have the form
$$
  \sum\limits_{s=1}^Nr_{is}\varphi^\mu_s=
  \sum\limits_{s=1}^N r_{is}\psi^\mu_s \quad (i=1,\dots,N;\, i\ne m). \eqno(18)
$$

Using the condition $\psi^\mu_{N+1}=\varphi^\mu_0=0,$ we rewrite
relations~(14) in the following form
$$
  \sum\limits_{s=1}^N r_{N+1,s}\psi^\mu_s=
  \sum\limits_{1\le i\le N+1,\,i\ne m+1}\gamma_{1i}
                        \sum\limits_{s=1}^N r_{i,s+1}\varphi^\mu_s. \eqno(19)
$$

The condition (13) implies that,
$$
  \sum\limits_{1\le i\le N+1,\,i\ne m+1}\gamma_{1i}
					\sum\limits_{s=1}^N r_{i,s+1}\varphi^\mu_s=
  \sum\limits_{s=1}^N r_{N+1,s}\varphi^\mu_s.
$$
And now, using (19), we obtain
$$
  \sum\limits_{s=1}^Nr_{N+1,s}\varphi^\mu_s=\sum\limits_{s=1}^N r_{N+1,s}\psi^\mu_s. 
								\eqno(20)
$$

Combining (18) and (20), we get the system of $N$ equations with $N$ unknowns
$$
  \sum\limits_{s=1}^Nr_{is}(\varphi^\mu_s-\psi^\mu_s)=0 \quad
  (i=1,\dots,N+1;\, i\ne m). \eqno(21)
$$

The rows of system (21) coincide with the linearly independent
rows $g_i\,(i=1,\dots,N+1;\, i\ne m).$ Hence $\varphi^\mu_s-\psi^\mu_s=0,$
i.e., $\varphi^\mu_s=\psi^\mu_s$ $(s=1,\dots,N;\, \mu=0,\dots,k-1).$
We have thus proved that $W^k_\gamma(0,N+1)\subset R_Q(\mathaccent23W^k(0,N+1)).$
$\Box$

\bigskip

{\bf Remark 1.} It can be given the following equivalent definition of the
subspace $W^k_\gamma(0,N+1).$
$W^k_\gamma(0,N+1)$ is the subspace of functions from $W^k(0,N+1)$
satisfying conditions
$$
  \begin{array}{l}
  u^{(\mu)}(0)=\sum\limits_{1\le i\le N+1,\,i\ne m'}
                    \gamma'_{1i}u^{(\mu)}(i),  \\
  u^{(\mu)}(m')=\sum\limits_{1\le i\le N,\,i\ne m'}
		      \gamma'_{2i}u^{(\mu)}(i),
  \end{array}
$$
where $m'$ is a fixed point from the set $\{1,\dots,N\},$
$\gamma'_{1i}(i=1,\dots,N+1,\, i\ne m'),$ 
$\gamma'_{2i}(i=1,\dots,N,\,i\ne m')$ are real numbers; 
$\mu=0,\dots,k-1;\; k\ge1.$

\bigskip

Let us introduce the sets
$$
  \begin{array}{rl}
      M&=\{u\in\mathaccent23W^1(0,N+1) : R_Qu\in W^2(0,N+1)\},\\
    M_k&=\{u\in\mathaccent23W^1(0,N+1) : u,\,R_Qu\in W^{k+2}(0,N+1)\}=\\
       &=\{u\in M : u,\,R_Qu\in W^{k+2}(0,N+1)\},
  \end{array}
$$
where $k=0,1,\dots$

These sets will play the role of the domains of the corresponding 
differential--difference operators.

We denote by $G^1_j\,(G^2_j)$ the $j$-th column of the 
$N\times(N+1)$-matrix obtained from $R_1$ by deleting the first (last) row
$(j=1,\dots,N+1).$ Notice that the conditions $\det R_1\ne 0,\, \det R_2=0$
imply that $G^1_1\ne 0,\,G^2_{N+1}\ne 0.$

The following lemma allows to find out the structure of the sets $M_k.$

\bigskip
{\bf Lemma 5.} {\sl For any $n\ge 2,$ we have:\newline
(a) Suppose that $G^1_1$ and $G^2_{N+1}$ are linearly independent. Then
$$
  \{v\in M : v,\,R_Qv\in W^n(0,N+1)\}=\mathaccent23W^n(0,N+1).
$$
(b) Suppose that $G^1_1$ and $G^2_{N+1}$ are linearly dependent. Then
$$
  \begin{array}{l}
    \{v\in M : v,\,R_Qv\in W^n(0,N+1)\} \\
    \qquad =\{v\in M : R_Qv\in W^n(0,N+1),\, (Uv)^{(\mu)}_{l+1}(0+0)=
							(Uv)^{(\mu)}_l(1-0)
    \quad (\mu=1,\dots,n-1)\},
  \end{array}
$$
where $l\in\{1,\dots,N\}$ is a point satisfying the following condition:
determinant of the matrix with the elements $r_{ij},$ where $1\le i,j\le N,\,
i\ne m,\,j\ne l,$ doesn't equal zero.
(By virtue of the linearly independence of the rows 
$g_i\, (i=1,\dots N,\,i\ne m),$ there really exists such a point $l$).
}    
\bigskip	
    	
{\it Proof.} First let us prove (a).

The inclusion $\mathaccent23W^n(0,N+1)\subset\{v\in M : v,\,R_Qv\in W^n(0,N+1)\}$ 
follows from Lemma 3. Let us prove the inverse inclusion.

Let $v\in\mathaccent23W^1(0,N+1)\cap W^n(0,N+1),\, R_Qv\in W^n(0,N+1).$ Then, using the 
notation of Theorem 1, for all $\mu=1,\dots,n-1,$ we obtain
$$
  \sum\limits_{s=1}^{N+1}r_{i+1,s}\varphi^\mu_{s-1}=
  \sum\limits_{s=1}^{N+1}r_{is}\psi^\mu_s \qquad (i=1,\dots,N).\eqno(22)
$$
Regrouping the summands in (22) and noticing that 
$r_{i+1,s+1}=r_{is}\,(1\le i,s\le N),$ we get
$$
  \sum\limits_{s=1}^N r_{is}(\varphi^\mu_s - \psi^\mu_s)=
  -r_{i+1,1}\varphi^\mu_0 + r_{i,N+1}\psi^\mu_{N+1} \qquad (i=1,\dots,N).
								   \eqno(23)
$$
Since $v\in W^n(0,N+1),$ we have $\varphi^\mu_s=\psi^\mu_s\, (s=1,\dots,N).$ 
Hence
$$
  -r_{i+1,1}\varphi^\mu_0+r_{i,N+1}\psi^\mu_{N+1}=0 \qquad (i=1,\dots,N).
$$
But the last relations are equivalent to the following
$$
  -G^1_1\varphi^\mu_0+G^2_{N+1}\psi^\mu_{N+1}=0.
$$
Thus, by virtue of the linearly independence of $G^1_1$ and $G^2_{N+1},$ we
have $\varphi^\mu_0=\psi^\mu_{N+1}=0.$ This implies that $v\in\mathaccent23W^n(0,N+1).$

Now let us prove (b). The inclusion 
$$
  \begin{array}{l}
    \{v\in M : v,\,R_Qv\in W^n(0,N+1)\} \\
    \subset\{v\in M : R_Qv\in W^n(0,N+1),\, (Uv)^{(\mu)}_{l+1}(0+0)=
							   (Uv)^{(\mu)}_l(1-0)
    \quad (\mu=1,\dots,n-1)\},
  \end{array}
$$
is obviously.

Let us prove the inverse inclusion. Let
$$
  v\in\{v\in M : R_Qv\in W^n(0,N+1),\, (Uv)^{(\mu)}_{l+1}(0+0)=
      (Uv)^{(\mu)}_l(1-0)\quad (\mu=1,\dots,n-1)\}.
$$

Note that it cannot be written ``$v^{(\mu)}(l-0)=v^{(\mu)}(l+0)$"
here and in the statement of the lemma because we don't know beforehand if the 
derivative of order $\mu$ for the function $v$ belongs to the
corresponding Sobolev space. Thus we have to write 
``$(Uv)^{(\mu)}_{l+1}(0+0)=(Uv)^{(\mu)}_l(1-0)$".

Since $G^1_1$ and $G^2_{N+1}$ are linearly dependent, there exist non--zero
real numbers $\alpha_1,\, \alpha_2$ such that 
$$
  \alpha_1 G^1_1+\alpha_2 G^2_{N+1}=0.\eqno(24)
$$
Now we shall show that, in this case,
$$
  \alpha_1 (Uv)^{(\mu)}_{N+1}(1-0)+\alpha_2 (Uv)^{(\mu)}_1(0+0)\equiv
  \alpha_1\psi^\mu_{N+1}+\alpha_2\varphi^\mu_0=0.\eqno(25)
$$

Denote $w=R_Qv.$ Since $(Uv)^{(\mu)}(t)=(R^{-1}_1Uw^{(\mu)})(t)\quad 
(t\in(0,1)),$ we can rewrite (25) in the form
$$
  \alpha_1\sum\limits_{i=1}^{N+1}\frac{B_{i,N+1}}{\det R_1}(Uw^{(\mu)})_i(1-0)+
  \alpha_2\sum\limits_{i=1}^{N+1}\frac{B_{i1}}{\det R_1}(Uw^{(\mu)})_i(0+0)=0.
 							\eqno(26)
$$
Since $B_{11}=B_{N+1,N+1}=\det R_2=0$, relation (26) has the form
$$
  \sum\limits_{i=1}^N (\alpha_1 B_{i,N+1}+\alpha_2 B_{i+1,1})w^{(\mu)}(i)=0. 
								      \eqno(27)
$$
Then, analyzing $B_{i,N+1},\;B_{i+1,1}$ and using (24), we see that
$\alpha_1B_{i,N+1}+\alpha_2B_{i+1,1}=0\; (i=1,\dots,N).$

Therefore (27) is identical, i.e., (25) is valid for any 
$$
  v\in\{v\in M : R_Qv\in W^n(0,N+1),\, (Uv)^{(\mu)}_{l+1}(0+0)=
                                     (Uv)^{(\mu)}_l(1-0)
                                   \quad (\mu=1,\dots,n-1)\}.
$$

Further, we have (likewise (23))
$$
  \sum\limits_{s=1}^N r_{is}(\varphi^\mu_s - \psi^\mu_s)=
  -r_{i+1,1}\varphi^\mu_0 + r_{i,N+1}\psi^\mu_{N+1} \qquad (i=1,\dots,N).
  								     \eqno(28)
$$
By virtue (24), (25), system (28) will have the form
$$
  \sum\limits_{s=1}^N r_{is}(\varphi^\mu_s - \psi^\mu_s)=0
  \qquad (i=1,\dots,N).\eqno(29)
$$
Since $\varphi^\mu_l=\psi^\mu_l$ and the $m$-th row of (29) is a linear
combination of the other ones, system (29) is equivalent to the following
$$
  \sum\limits_{1\le s\le N,\,s\ne l}
       r_{is}(\varphi^\mu_s - \psi^\mu_s)=0
  \qquad (i=1,\dots,N;\,i\ne m).\eqno(30)
$$
Thus we have the system of $(N-1)$ equations with $(N-1)$ unknowns. Selection
of point $l$ implies non-singularity of the matrix of system (30). This system
has a unique trivial solution. Hence, for any $\mu=0,\dots,n-1,$ we get 
$\varphi^\mu_s=\psi^\mu_s\, (s=1,\dots\ N).$ Therefore $v\in W^n(0,N+1)$ and
thus Lemma~5 is proved. $\Box$
\bigskip

Let $R_Q^k:W^{k+2}(0,N+1)\to W^{k+2}(0,N+1)$ be a bounded operator defined by
${\cal D}(R_Q^k)=M_k,\quad R_Q^kv=R_Qv\, (v\in {\cal D}(R_Q^k)),$ 
where $k\ge 0.$

\bigskip
{\bf Theorem 2.} {\sl The operator $R_Q^k\, (k\ge 0)$ is Fredholm, 
$\dim\ker(R_Q^k)=0,$\newline
$
  \mbox{\rm codim\,}\mbox{\rm Im\,}(R^k_Q)=\left\{
  \begin{array}{ll} 2(k+2), &\mbox {if $G^1_1,\, G^2_{N+1}$ are linearly 
								independent,}\\
                   k+3, &\mbox {if $G^1_1,\, G^2_{N+1}$ are linearly dependent.}
  \end{array}\right. 
$
}
\bigskip

{\it Proof.} Let $G^1_1,\,G^2_{N+1}$ be linearly independent. In this case,
by virtue of Lemma 5, the domain $M_k$ of the operator $R_Q^k$ coincides 
with the space $\mathaccent23W^{k+2}(0,N+1).$ By virtue of Theorem 1, 
the operator $R_Q^k$ maps $M_k$ onto $W^{k+2}_\gamma(0,N+1)$ in a one--to--one
manner. This implies that $\dim\ker(R_Q^k)=0.$

Now let us find $\mbox{\rm codim\,}\mbox{\rm Im\,}(R_Q^k).$ We consider the equation
$$
  R_Q^ku=w\qquad (w\in W^{k+2}(0,N+1)).\eqno(31)
$$ 

Theorem 1 implies that equation (31) has a solution $u\in M_k=\mathaccent23W^{k+2}(0,N+1)$
iff $w\in W^{k+2}_\gamma(0,N+1),$ i.e., iff $w$ satisfies the conditions
$$
  \begin{array}{l}
  w^{(\mu)}(N+1)=\sum\limits_{1\le i\le N+1,\,i\ne m+1}\gamma_{1i}w^{(\mu)}(i-1), \\
  w^{(\mu)}(m)=\sum\limits_{1\le i\le N,\,i\ne m}\gamma_{2i}w^{(\mu)}(i) 
						\qquad (\mu=0,\dots,k+1).
  \end{array}
$$

We introduce $2(k+2)$ linear functionals $F_{j\mu}\, (j=0,1;\,\mu=0,\dots,k+1)$
by the formulas
$$
  \begin{array}{l}
  F_{0\mu}(w)=w^{(\mu)}(N+1)-\sum\limits_{1\le i\le N+1,\,i\ne m+1}
						\gamma_{1i}w^{(\mu)}(i-1), \\
  F_{1\mu}(w)=w^{(\mu)}(m)-\sum\limits_{1\le i\le N,\,i\ne m}\gamma_{2i}w^{(\mu)}(i).
  \end{array}\eqno(32)
$$
By virtue of the trace theorem (for example, see [2]), $F_{j\mu}$ are 
continuous functionals over $W^{k+2}(0,N+1)$. It is not hard to check that 
$F_{j\mu}$ are linearly independent.

From the Riesz theorem it follows that $F_{j\mu}(w)=(w,f_{j\mu})_{W^{k+2}(0,N+1)},$
where $f_{j\mu}\in W^{k+2}(0,N+1)\, (j=0,1;\, \mu=0,\dots,k+1)$ are linearly 
independent functions. This implies that $\mbox{\rm codim\,}\mbox{\rm Im\,}(R_Q^k)=2(k+2).$

Now we consider the other case. Let $G_1^1,\,G^2_{N+1}$ be linearly dependent.
Since ${\cal D}(R_Q^k)\subset\mathaccent23W^1(0,N+1),$ $R_Q$ maps
$\mathaccent23W^1(0,N+1)$ onto $W^1_\gamma(0,N+1)$ in a one--to--one manner, and
$R_Q\supset R_Q^k,$ it follows that $\dim\ker(R_Q^k)=0.$

Let us find $\mbox{\rm codim\,}\mbox{\rm Im\,}(R_Q^k).$ We consider the equation
$$
  R_Q^kv=w\qquad (w\in W^{k+2}(0,N+1)).\eqno(33)
$$ 

From Theorem~1 and Lemma~5, it follows that equation~(33) has a solution 
$v\in M_k$ iff $w$ satisfies the conditions
$$
  w(N+1)=\sum\limits_{1\le i\le N+1,\,i\ne m+1}\gamma_{1i}w(i-1),\quad 
  w(m)=\sum\limits_{1\le i\le N,\,i\ne m}\gamma_{2i}w(i), \eqno(34)
$$
$$
    (Uv)^{(\mu)}_{l+1}(0+0)=(Uv)^{(\mu)}_l(1-0)\qquad (\mu=1,\dots,k+1).
						 	\eqno(35)
$$
Since
$$
  \begin{array}{l}
   (Uv)^{(\mu)}_{l+1}(0+0)=\sum\limits_{i=1}^{N+1}\frac{B_{i,l+1}}{\det R_1}
				       (Uw^{(\mu)}_i(0+0)),\\
   (Uv)^{(\mu)}_l(1-0)=\sum\limits_{i=1}^{N+1}\frac{B_{il}}{\det R_1}
				       (Uw^{(\mu)}_i(1-0)),
  \end{array}
$$
conditions~(35) will have the form
$$
 \sum\limits_{i=1}^{N+1}B_{i,l+1}w^{(\mu)}(i-1)=
 \sum\limits_{i=1}^{N+1}B_{il}w^{(\mu)}(i).
$$
And, after regrouping of the summands, we obtain
$$
   B_{1,l+1}w^{(\mu)}(0)+\sum\limits_{i=1}^{N+1}(B_{i+1,l+1}-B_{il})w^{(\mu)}(i)-
   B_{N+1,l}w^{(\mu)}(N+1)=0 \quad  (\mu=1,\dots,k+1). \eqno(36)
$$
Thus a solution $u$ of equation (33) belongs to $M_k$ iff $w$ satisfies 
conditions~(34) and (36).
Further, as above, we can introduce $k+3$ linear continuous functionals over
$W^{k+2}(0,N+1)$, corresponding conditions (34), (36), and prove that they are
linearly independent. (To prove it we need the condition $B_{N+1,l}\ne 0$
which follows from the conditions on the point $l$.) And, as above, using the 
Riesz theorem, we get
$\mbox{\rm codim\,}\mbox{\rm Im\,}(R_Q^k)=k+3.$ $\Box$
\bigskip          
\begin{center}
\large\bf \S2. The boundary value problem for the differential--difference\\
equation with homogeneous boundary conditions
\end{center}

We consider the {\it differential--difference equation}
$$
  -(Rv)''(t)+(A_1v)(t)=f_0(t)\qquad (t\in(0,N+1))\eqno(37)
$$
with homogeneous boundary conditions
$$ 
  v(t)=0\qquad (t\in [-N,0]\cup [N+1,2N+1]).\eqno(38)
$$
Here $R:L_2({\bf R})\to L_2({\bf R})$ is the difference operator defined by
$$
  (Rv)(t)=\sum\limits_{j=-N}^Nb_jv(t+j),
$$
$b_j\in{\bf R};\,N\in{\bf N};\, A_1:\mathaccent23W^1(0,N+1)\to L_2(0,N+1)$ is a linear 
bounded operator; $f_0\in L_2(0,N+1).$ One can easily reduce a 
differential--difference equation with non-homogeneous boundary conditions
to differential--difference equation with homogeneous boundary conditions
(see \S3). Therefore, without loss of generality, we can study the equation~(37)
with homogeneous boundary conditions~(38). 

Since the shifts $t\to t+j$
can map the points of $(0,N+1)$ into the set $[-N,0]\cup[N+1,2N+1],$ we consider
the boundary conditions for the equation~(37) not only at the ends of the
interval $(0,N+1),$ but also on the set $[-N,0]\cup[N+1,2N+1].$

Let $A_R:L_2(0,N+1)\to L_2(0,N+1)$ be the unbounded operator given by
$$
  \begin{array}{rl}
    {\cal D}(A_R)&=M=\{v\in\mathaccent23W^1(0,N+1) : R_Qv\in W^2(0,N+1)\},\\
    A_Rv&=-(R_Qv)''(t)+A_1v \qquad (v\in{\cal D}(A_R)).
  \end{array}
$$

{\bf Definition 1.} A function $v\in{\cal D}(A_R)$ is called
{\it a generalized solution} for problem (37), (38) if 
$
  A_Rv=f_0.
$ 
\bigskip

{\bf Theorem 3.} {\sl The operator $A_R$ is Fredholm and $\mbox{ind\,} A_R=0.$}
\bigskip

To prove Theorem 3 we first consider the bounded operator 
$A:W^2(0,N+1)\cap W^1_\gamma(0,N+1)\to L_2(0,N+1)$ defined by the formula
$$
  Au=-u''+A_1R_Q^{-1}u.
$$
Here we suppose that the space $W^2(0,N+1)\cap W^1_\gamma(0,N+1)$ 
has a topology of the space $W^2(0,N+1)$.
Let us prove the following lemma.

\bigskip
{\bf Lemma 6.} {\sl The bounded operator $A$ is Fredholm and $\mbox{ind\,} A=0.$}
\bigskip

{\it Proof.} We introduce the bounded operator 
$A_2:W^2(0,N+1)\cap W^1_\gamma(0,N+1)\to L_2(0,N+1)$
defined by the formula
$$
  A_2u=u''(t).
$$
Here we also suppose that the space 
$W^2(0,N+1)\cap W^1_\gamma(0,N+1)$ has a topology of $W^2(0,N+1)$.

Thus we have $A=-A_2+A_1R_Q^{-1}.$
We show that the operator $A_2$ is Fredholm and $\mbox{ind\,} A_2=0.$

It is clear that the homogeneous equation $A_2u\equiv u''(t)=0$ has a
class of solutions $u(t)=c_1t+c_2$ from $W^2(0,N+1).$ Therefore $u$ belongs to 
$\ker(A_2)$ iff $u$ satisfies conditions (5), (6) (for $\mu=0$)
$$
  \begin{array}{l}
  c_1[N+1-\sum\limits_{2\le i\le N+1,\,i\ne m+1}\gamma_{1i}(i-1)]+
   c_2[1-\sum\limits_{1\le i\le N+1,\,i\ne m+1}\gamma_{1i}]=0,\\
  c_1[m-\sum\limits_{1\le i\le N,\,i\ne m}\gamma_{2i}i]+
   c_2[1-\sum\limits_{1\le i\le N,\,i\ne m}\gamma_{2i}]=0.
  \end{array}\eqno(39)
$$
Parallel with the homogeneous equation, we shall consider the non-homogeneous
equation
$$
  A_2v\equiv v''(t)=f(t)\qquad (f\in L_2(0,N+1)).
$$
For any function $f\in L_2(0,N+1)$, there exists a class of solutions
$v(t)=d_1t+d_2+\int\limits_0^t (t-\tau)f(\tau)\,d\tau$ from $W^2(0,N+1).$
Therefore $v$ belongs to the domain of the operator $A_2$ iff $v$ satisfies 
conditions~(5), (6) (for $\mu=0$)
$$
  \begin{array}{l}
    d_1[N+1-\sum\limits_{2\le i\le N+1,\,i\ne m+1}\gamma_{1i}(i-1)]+
    d_2[1-\sum\limits_{1\le i\le N+1,\,i\ne m+1} \gamma_{1i}]                                 \\
    \qquad=\sum\limits_{1\le i\le N,\,i\ne m}\gamma_{1,i+1}
          \int\limits_0^i (i-\tau)f(\tau)\,d\tau - 
          \int\limits_0^{N+1}(N+1-\tau)f(\tau)\,d\tau,               \\
    d_1[m-\sum\limits_{1\le i\le N,\,i\ne m} \gamma_{2i}i]+
    d_2[1-\sum\limits_{1\le i\le N,\,i\ne m} \gamma_{2i}]=
     \sum\limits_{1\le i\le N,\,i\ne m}\gamma_{2i}
          \int\limits_0^i (i-\tau)f(\tau)\,d\tau - 
          \int\limits_0^m(m-\tau)f(\tau)\,d\tau.
  \end{array}\eqno(40)
$$

It is clear that $\Phi_i(f)=(f,\phi_i)_{L_2(0,N+1)}\equiv
\int\limits_0^i (i-\tau)f(\tau)\,d\tau\; (i=1,\dots,N+1)$ are the linear 
continuous functionals over $L_2(0,N+1)$ (here $\phi_i(\tau)=(i-\tau)I(i-\tau),$
where $I(t)=1,\, t\ge 0;\; I(t)=0,\, t<0).$

It is not hard to prove that the functionals $\Phi_i\, (i=1,\dots,N+1)$ are
linearly independent. This implies that
$$
  \begin{array}{l}
   F_1(f)=\sum\limits_{1\le i\le N,\,i\ne m}\gamma_{1,i+1}\Phi_i(f)-\Phi_{N+1}(f),\\
   F_2(f)=\sum\limits_{1\le i\le N,\,i\ne m}\gamma_{2i}\Phi_i(f)-\Phi_m(f)
  \end{array}
$$
are non-zero linearly independent continuous functionals over $L_2(0,N+1).$

Thus system (40) will have the form
$$
  \begin{array}{l}
    d_1[N+1-\sum\limits_{1\le i\le N+1,\,i\ne m+1} \gamma_{1i}(i-1)]+
    d_2[1-\sum\limits_{1\le i\le N+1,\,i\ne m+1} \gamma_{1i}]=F_1(f),\\                            \\
    d_1[m-\sum\limits_{1\le i\le N,\,i\ne m} \gamma_{2i}i]+
    d_2[1-\sum\limits_{1\le i\le N,\,i\ne m} \gamma_{2i}]=F_2(f).  
  \end{array}\eqno(41)
$$

We analyse system (39) and system (41) simultaneously. Notice that the
matrix of system (39) coincides with the matrix of system (41). Denote
this matrix by $\cal M.$ Let us consider three cases.

1. $\mbox{\rm Rank\,}({\cal M})=2.$ It is easy to see that we have 
$\dim\ker(A_2)=0,\, \mbox{\rm codim\,}\mbox{\rm Im\,}(A_2)=0,$ i.e., $\mbox{ind\,} A_2=0.$

2. $\mbox{\rm Rank\,}({\cal M})=1.$ Clearly, $\dim\ker(A_2)=1.$ Using the Riesz 
theorem, we obtain $\mbox{\rm codim\,}\mbox{\rm Im\,}(A_2)=1.$ Hence, in this case, we also have 
$\mbox{ind\,} A_2=0.$

3. $\mbox{\rm Rank\,}({\cal M})=0.$ In this case, we see that $\dim\ker(A_2)=2.$ Using again
the Riesz theorem, we obtain $\mbox{\rm codim\,}\mbox{\rm Im\,}(A_2)=2,$ i.e., $\mbox{ind\,} A_2=0.$

Thus we have proved that $A_2$ is Fredholm and $\mbox{ind\,} A_2=0.$

It is not hard to check that the operator 
$A_1R_Q^{-1}:W^2(0,N+1)\cap W^1_\gamma(0,N+1)\to L_2(0,N+1)$ is bounded if the
space $W^2(0,N+1)\cap W^1_\gamma(0,N+1)$ has a topology of $W^1(0,N+1)$.
Therefore, by virtue of the compactness of the embedding operator from 
$W^2(0,N+1)$ into $W^1(0,N+1)$, the operator 
$A_1R_Q^{-1}:W^2(0,N+1)\cap W^1_\gamma(0,N+1)\to L_2(0,N+1)$ 
is compact if the space $W^2(0,N+1)\cap W^1_\gamma(0,N+1)$ has a topology of 
$W^2(0,N+1)$. Using the theorem about the compact perturbations (see [1], theorem 16.4), 
we have that the operator $A=-A_2+A_1R_Q^{-1}$ is Fredholm and $\mbox{ind\,} A=0.$ $\Box$
\bigskip

Now let us prove Theorem 3.

{\it Proof of Theorem 3.} The operator $A_R$ can be presented as a composition 
$A_R=A\mathaccent "7E R_Q,$ where\newline
$A:W^2(0,N+1)\cap W^1_\gamma(0,N+1)\to L_2(0,N+1)$ is given by
$$
  Au=-u''+A_1R_Q^{-1}u,
$$
$\mathaccent "7E R_Q:L_2(0,N+1)\to W^2(0,N+1)\cap W^1_\gamma(0,N+1)\subset 
W^2(0,N+1)$ is given by
$$
  \begin{array}{rl}
      {\cal D}(\mathaccent "7E R_Q)&={\cal D}(A_R)=M,\\
      \mathaccent "7E R_Qu&=R_Qu \qquad (u\in{\cal D}(\mathaccent "7E R_Q)).
  \end{array}
$$
By virtue of Lemma 6 and Theorem 1, the operators $A$ and $\mathaccent "7E R_Q$
are Fredholm and $\mbox{ind\,} A=\mbox{ind\,}\mathaccent "7E R_Q=0.$ Hence the operator
$A_R=A\mathaccent "7E R_Q$ is also Fredholm and $\mbox{ind\,} A_R=0$ (see [1], 
theorem 12.2). $\Box$
\begin{center}\large\bf
\S3. Smoothness of generalized solutions to boundary value problem
\end{center}

It is known that the smoothness of generalized solutions of 
differential--difference equations can be broken even for infinitely
differentiable right hand sides of equations. But there exists the following
result.

\bigskip
{\bf Theorem 4.} {\sl Let $f_0\in W^k(0,N+1)$ and $v$ be a generalized solution
for boundary value problem~(37), (38) such that $A_1v\in W^k(0,N+1)$.

Then $v\in W^{k+2}(Q_s),\; s=1,\dots,N+1.$ }
\bigskip
		 
{\it Proof.} The proof follows from Lemma 4.
\bigskip

To obtain a smoothness of generalized solutions it is necessary to impose some
additional conditions on right hand side of the equation (and on the boundary 
functions, in the case of non-homogeneous boundary conditions). 
Now we shall find out a type of these conditions for the case of 
the homogeneous boundary value problem.

We consider the bounded operator $A_R^k:W^{k+2}(0,N+1)\to W^k(0,N+1)$ given by
$$
  \begin{array}{l}
    {\cal D}(A^k_R)=M_k,\\
    A_R^kv=-(R_Qv)''(t),
  \end{array}
$$ 
and the bounded operator $B_R^k:\mathaccent23W^{k+2}(0,N+1)\to W^k(0,N+1)$ defined by the
formula $B_R^kv=-(R_Qv)''(t).$

Note that, by virtue of Lemma 5, $A_R^k$ coincides with $B_R^k$ if 
$G_1^1,\;G_{N+1}^2$ are linearly independent.

\bigskip
{\bf Theorem 5.} {\sl The operator $A_R^k\; (k\ge 0)$ is Fredholm, 
$\dim\ker(A_R^k)=0,$\newline 
$\mbox{\rm codim\,}\mbox{\rm Im\,}(A_R^k)=\left\{
  \begin{array}{ll} 2(k+1), & \mbox{if $G_1^1,\,G_{N+1}^2$ are linearly 
								independent,}\\
                   k+1, & \mbox{if $G_1^1,\,G_{N+1}^2$ are linearly dependent.}
  \end{array}\right. $ }
\bigskip

{\it Proof.} First we prove that $\dim\ker(A_R^k)=0.$
Let $v\in\ker(A^k_R).$ Then $(R_Qv)''(t)=0.$ Hence 
$(R_Qv)(t)=c_1+c_2t.$ Since $\det R_1\ne0$, we obtain
$$
  v(t)=U^{-1}R_1^{-1}U(c_1+c_2t)\qquad (t\in(0,N+1)).
$$
Thus a function $v$ is piecewise linear on the interval $(0,N+1).$ Therefore
$v\in W^2(0,N+1)\cap\mathaccent23W^1(0,N+1)$ if and only if $v(t)=0,$ i.e., 
$\dim\ker(A_R^k)=0.$

Let us present the operator $A_R^k$ as a composition
$A_R^k=A_2R_Q^k.$ Here $R_Q^k: W^{k+2}(0,N+1)\to W^{k+2}(0,N+1)$ is the operator 
introduced in \S1, $A_2:W^{k+2}(0,N+1)\to W^k(0,N+1)$ is the bounded
operator defined by the formula $(A_2v)(t)=-v''(t).$ It is obvious that
$A_2$ is Fredholm and $\mbox{\rm ind\,} A_2=2.$ Therefore, using Theorem~2 and the theorem
about a composition of Fredholmian operators (see [1], theorem 12.2), we
obtain the statement of Theorem 5.
\bigskip

{\bf Theorem 6.} {\sl The operator $B_R^k\; (k\ge 0)$ is Fredholm, 
$\dim\ker(B_R^k)=0,\; 
\mbox{\rm codim\,}\mbox{\rm Im\,}(B_R^k)=2(k+1).$ }
\bigskip

{\it Proof.} The idea of the proof is analogous to the previous proof.
\bigskip

Now we shall generalize these results to
the case of the boundary value problem with non-homogeneous boundary conditions.

We consider the differential--difference equation
$$
  -(Ry)''(t)+A_1y=f_0(t)\qquad (t\in(0,N+1))\eqno(42)
$$
with non-homogeneous boundary conditions
$$
  \left\{
  \begin{array}{ll}
     y(t)=f_1(t) & (t\in [-N,0]),\\
     y(t)=f_2(t) & (t\in [N+1,2N+1]),	
  \end{array}\right.\eqno(43)
$$
where
$$
  (Ry)(t)=\sum\limits_{j=-N}^N b_jy(t+j),
$$
$b_j\in{\bf R},$ $N$ is a natural number; $A_1:W^1(-N,2N+1)\to L_2(0,N+1)$ is a
linear bounded operator, $f=(f_0,f_1,f_2)\in{\cal W}(-N,2N+1)=L_2(0,N+1)\times
W^1(-N,0)\times W^1(N+1,2N+1).$

We introduce the linear unbounded operator ${\cal L}:L_2(-N,2N+1)\to
{\cal W}(-N,2N+1)$ with the domain ${\cal D}({\cal L})=\{y\in W^1(-N,2N+1) : 
P_QRy\in W^2(0,N+1)\}$ by the formula
$$
  {\cal L}y=\left(-(P_QRy)''+A_1y,y|_{(-N,0)},y|_{(N+1,2N+1)}\right).
$$
{\bf Definition 2.} A function $y\in{\cal D}({\cal L})$ is called {\it
a generalized solution} for problem (42), (43) if ${\cal L}y=(f_0,f_1,f_2).$
\bigskip

To obtain the smoothness of the generalized solution in the interval
$(-N,\,2N+1)$ we suppose that $A_1:W^{k+1}(-N,\,2N+1)\to W^k(0,N+1)$ is a
bounded operator and
$ f=(f_0,f_1,f_2)\in{\cal W}^k(-N,\,2N+1)=W^k(0,N+1)\times W^{k+2}(-N,0)\times
W^{k+2}(N+1,\,2N+1).$

We consider the linear bounded operator ${\cal L}_B:W^{k+2}(-N,\,2N+1)\to
{\cal W}^k(-N,\,2N+1)$ by the formula 
$$
  {\cal L}_By={\cal L}y\qquad (y\in W^{k+2}(-N,\,2N+1)).
$$

{\bf Theorem 7.} {\sl The operator ${\cal L}_B$ is Fredholm and 
$\mbox{\rm ind\,} {\cal L}_B=-2(k+1).$ }
\bigskip

{\it Proof.} By virtue of the compactness of the imbedding operator from
$W^{k+2}(-N,\,2N+1)$ into $W^{k+1}(-N,\,2N+1),$ the operator 
$A_1:W^{k+2}(-N,\,2N+1)\to W^k(0,N+1)$ is compact. Therefore, by theorem 16.4,
[1], it suffices to prove Theorem 7 in the case $A_1=0.$

Let us assume now that $A_1=0$

We introduce the function
$$
\psi(t)=\left\{
        \begin{array}{l}
         f_1(t) \qquad (t\in [-N,0]),\\
	 f_2(t) \qquad (t\in [N,\,2N+1]),\\
         \eta(t)\sum\limits_{i=0}^{k+1}f_1^{(i)}(0)t^i/i!+\eta(t-N-1)\sum\limits_{i=0}^{k+1}
      			         f_2^{(i)}(N+1)(t-N-1)^i/i!\\
         \phantom{aaaaaaaaaaaaaaaaaaaaaaaaaaaaaaaaaaaaaaaaa}
							(t\in(0,N+1)),
        \end{array}\right.
$$
where $\eta\in \dot C^\infty({\bf R}),\; \eta(t)=1\; (|t|<1/4),\;\eta(t)=0\;
(|t|>1/3).$ It is clear that $\psi\in W^{k+2}(-N,\,2N+1).$ Denote 
$w=y-\psi\in W^{k+2}(0,N+1)$ ($y\in W^{k+2}(-N,\,2N+1)$). We see that the
equation ${\cal L}_By=f\;(f\in{\cal W}^k(-N,\,2N+1))$ has a solution 
$y\in W^{k+2}(-N,\,2N+1)$ iff $w$ belongs to $\mathaccent23W^{k+2}(0,N+1)$ and is a 
solution of the equation
$$
  B_R^kw=f_0+(R\psi)''.\eqno(44)
$$
By Theorem 6, equation (44) has a solution if and only if
$$
  \left(f_0+(R\psi)'',\varphi_j\right)_{W^k(0,N+1)}=0\qquad 
	\left(j=1,\dots,2(k+1)\right),\eqno(45)
$$
where $\varphi_j\in W^k(0,N+1)$ are linearly independent functions.

From the trace theorem and the Riesz theorem it follows that conditions~(45)
will have the form
$$
  (f,G_j)_{{\cal W}^k(-N,\,2N+1)}=0\qquad \left(j=1,\dots,2(k+1)\right),
								     \eqno(46)
$$
where $f=(f_0,f_1,f_2),$ vector-valued functions $G_j=(\varphi_j,B_1\varphi_j,
B_2\varphi_j)$ are linearly independent (here $B_1:W^k(0,N+1)\to W^{k+2}(-N,0),\;
B_2:W^k(0,N+1)\to W^{k+2}(N+1,\,2N+1)$ are linear bounded operators). Thus for
$A_1=0$ the equation ${\cal L}_By=f$ has a solution $y\in W^{k+2}(-N,\,2N+1)$
for $f\in{\cal W}^k(-N,\,2N+1)$ if and only if conditions~(46) are fulfilled.

Furthermore, by Theorem 6, $\dim\ker({\cal L}_B)=0.$ $\Box$
\bigskip

If we demand the smoothness of the solution only in the interval $(0,N+1)$,
we can weaken the conditions of orthogonality in some cases. 

To formalize this statement we suppose that $A_1:W^1(-N,\,2N+1)\to W^k(0,N+1)$ 
is a compact operator. Let us introduce the unbounded operator 
${\cal L}_A:W^1(-N,2N+1)\to W^k(0,N+1)\times W^1(-N,0)\times W^1(N+1,\,2N+1)$ 
with the domain 
${\cal D}({\cal L}_A)=\{y\in W^1(-N,\,2N+1):P_Qy,P_QRy\in W^{k+2}(0,N+1)\}$ 
by the formula 
$$
  {\cal L}_Ay={\cal L}y\qquad (y\in{\cal D}({\cal L}_A)).
$$

{\bf Theorem 8.} {\sl The operator 
${\cal L}_A$ is Fredholm and \newline
  $\mbox{\rm ind\,} {\cal L}_A=\left\{
  \begin{array}{ll} -2(k+1), & \mbox{if $G_1^1,\,G_{N+1}^2$ are linearly 
								independent,}\\
         -(k+1),  & \mbox{if $G_1^1,\,G_{N+1}^2$ are linearly dependent.}
  \end{array}\right.$ }
\bigskip

{\it Proof.} The proof is analogous to the previous proof. The main distinction
refers to the operator on left hand side of equation~(44) to which we reduce 
boundary value problem~(42), (43).

In this case, $A_R^k$ takes the place of $B_R^k.$ $\Box$
\bigskip

{\bf Remark 2.} Using Lemma 5, we can show that ${\cal D}({\cal L}_A)=
W^{k+2}(-N,\,2N+1)$ if $G_1^1,\,G_{N+1}^2$ are linearly independent. In this
case, a generalized solution has a proper smoothness in the whole interval
$(-N,\,2N+1).$
\bigskip

Thus we see that the smoothness of generalized solutions of the boundary value
problem to differential--difference equations is not broken in the interval 
$(0,N+1)$ (in the interval $(-N,\,2N+1)$) if we impose not only the conditions 
of smoothness but also some conditions of orthogonality on the right hand side 
of the differential--difference equation and on the boundary functions.
\begin{center}\large\bf
Bibliography
\end{center}

[1] S.\,G.\,Kre\v{\i}n, {\it Linear Equations in Banach Spaces}, Nauka, Moscow, 1971
(Russian); English translation: Birkh\"auser, Boston, 1982.

[2] V.\,P.\,Mikhailov, {\it Partial Differential Equations}, Nauka, Moscow, 1983
(Russian).

[3] A.\,L.\,Skubachevski\v{\i}, {\it Elliptic Functional Differential Equations 
and Applications}, Basel--Boston--Berlin, Birkh\"auser, 1997. 
\end{document}